\documentclass[a4paper, 11pt, reqno]{amsart}

\usepackage[usenames,dvipsnames]{color}
\usepackage{amssymb, amsmath, latexsym, graphics, graphicx, pstricks}

\newcommand\GreenL{\mathrel{\mathcal L}}
\newcommand\GreenR{\mathrel{\mathcal R}}

\newcommand\GreenD{\mathrel{\mathcal D}}
\newcommand\GreenJ{\mathrel{\mathcal J}}

\newcommand\ft{\mathbb{FT}}

\newcommand\transpose{\mathrm T}

\renewcommand\ker{\operatorname{Ker}}

\DeclareMathOperator\row{Row}
\DeclareMathOperator\col{Col}
\DeclareMathOperator\mx{M}
\DeclareMathOperator\ann{Ann}

\newtheorem{theorem}{Theorem}[section]

\newtheorem{proposition}[theorem]{Proposition}

\newtheorem{corollary}[theorem]{Corollary}
\newtheorem{example}[theorem]{Example}
\newtheorem{definition}[theorem]{Definition}

\begin{document}
	\title[Exact rings and semirings]{Exact rings and semirings}

	\keywords{rings; semirings; Hahn-Banach separation; tropical semiring; Green's relations}
	\subjclass[2010]{14T05, 16Y60, 20M10}

	\maketitle

	\begin{center}
		DAVID WILDING\footnote{Email \texttt{david.wilding-2@postgrad.manchester.ac.uk}.
		Research supported by a University of Manchester Faculty of Engineering and Physical Sciences Dean's Award.},
		MARIANNE JOHNSON\footnote{Email \texttt{Marianne.Johnson@maths.manchester.ac.uk}.
		Research supported by EPSRC grant number EP/H000801/1 (\textit{Multiplicative Structure of Tropical Matrix Algebra}).}
		and MARK KAMBITES\footnote{Email \texttt{Mark.Kambites@manchester.ac.uk}.
		Research also supported by EPSRC grant number EP/H000801/1.}

		\medskip

		School of Mathematics, \ University of Manchester, \\
		Manchester M13 9PL, \ England.

		\keywords{}
		\thanks{}
	\end{center}

	\numberwithin{equation}{section}

	\begin{abstract}
		We introduce and study an abstract class of semirings, which we call
		exact semirings, defined by a
		Hahn-Banach-type separation property on modules. Our motivation comes
		from the tropical semiring, and in particular a desire to understand the often
		surprising extent to which it behaves like a field. The definition of
		exactness abstracts an elementary property common to both fields and the tropical
		semiring, which we believe is fundamental to explaining this similarity.
		The class of exact semirings turns out to include many other important examples of both
		rings (proper quotients of principal ideal domains, matrix rings and finite
		group rings over these and over fields), and semirings (the Boolean
		semiring, generalisations of the tropical semiring, matrix semirings and
		group semirings over these).
	\end{abstract}

	\section{Introduction}
	A \textit{semiring} is an algebraic structure satisfying the usual axioms
	for a (not necessarily commutative) ring, but without the requirement that
	addition be invertible. Aside from rings, well-studied examples include the
	\textit{Boolean semiring} and the \textit{tropical semiring}.
	The latter is the algebraic structure formed by the real numbers (sometimes
	considered with $-\infty$ adjoined) under the operations of addition and
	maximum, which play the roles of semiring multiplication and (non-invertible)
	semiring addition respectively. It has a huge array of applications in areas including
	enumerative algebraic geometry and discrete event systems, and has
	been independently rediscovered many times by people working in these areas.

	A lack of additive inverses is clearly a radical departure from the
	definition of a field, and \textit{a priori} one would not expect linear algebra
	over any non-ring semiring (even a semifield) to behave much like classical
	linear algebra over a field. For example, as a direct consequence of the
	non-invertibility of addition, invertible tropical matrices are
	extremely sparse \cite[Theorem 1.1.3]{Butkovic10} and have little do with the structure of matrices in
	general. There is also no single well-behaved notion of
	\textit{rank} for tropical matrices, partly as a consequence of which
	even submodules of free modules can be exceedingly complex objects.
	However, as one delves deeper into the theory, one is often
	struck by surprising large-scale similarities between the tropical and
	the field case.

	A case in point is Green's $\GreenD$ relation (a key tool for
	capturing the structure of left and right ideals in semigroup theory
	\cite{Green51,Howie95}) for matrices.
	It is
	folklore (\cite[Theorem 1.4]{Putcha88} and \cite[Lemma 2.1]{Okninski98}) that two matrices (square and of
	the same size) over a field
	are $\GreenD$-related if and only if they have the same rank.
	Since the isomorphism type of a vector space is uniquely determined by its rank,
	another way of saying this is that they are $\GreenD$-related if and only
	if they have isomorphic column (or equivalently, row) spaces. In
	\cite{K_tropd}, Hollings and the third author showed that two matrices
	over the tropical semiring are also $\GreenD$-related if and only if they
	have (tropical linearly) isomorphic column spaces.

	A key ingredient in \cite{K_tropd} was a lemma establishing a kind of
	elementary Hahn-Banach separation property, which we shall here call
	\textit{exactness}. This property, which has a number of
	equivalent formulations (see Section~\ref{sec_exact} below), arises in
	the tropical case from the phenomenon of \textit{tropical matrix duality}
	\cite{Cohen04,Develin04,K_tropd}, and also holds (for completely different
	and much more elementary reasons) in fields. It has a number of interesting
	consequences, and indeed we believe it may explain much of the apparent
	commonality between these ostensibly quite different structures.

	From an algebraic perspective, one is drawn to ask whether other
	important semirings (including rings) are exact in this sense, and hence whether they can
	be expected to share that behaviour which is common to fields and the
	tropical semiring. If so, it would be natural to develop an
	abstract algebraic theory of such semirings. In fact the class of exact
	semirings does indeed turn out to encompass a number
	of important semirings. These include the Boolean semiring, proper quotients
	of principal ideal domains (such as $\mathbb Z/n\mathbb Z$ for all non-zero $n\in\mathbb Z$)
	and, more generally, self-injective rings.
	It is also includes matrix semirings over exact semirings, and finite (and
	in some cases infinite) group semirings over exact semirings. On the
	other hand, many important rings do \textit{not}
	have this property: indeed we shall see that an integral domain is exact
	only if it is a field, so $\mathbb{Z}$, for example, is \textit{not} exact.

	In addition to this introduction, this article comprises seven
	sections. Section~\ref{sec_prelim} recaps the definitions of semirings
	and various associated concepts. Section~\ref{sec_exact} introduces
	exact semirings, proves the equivalence of a number of different
	characterizations, and establishes some basic properties.
	Section~\ref{sec_rings} considers what rings have the exactness property,
	showing that a proper quotient of a principal ideal domain is always
	exact, but an integral domain cannot be exact unless it is a field.
	Section~\ref{sec_matrix} shows
	that the class of exact semirings is closed under taking matrix semirings;
	as a consequence we also deduce that a finite
	group semiring over an exact semiring is exact.
	Section~\ref{sec_anti} considers a class of semirings which are in
	a certain sense \textit{anti-isomorphic} to themselves, showing that these
	semirings are also exact: these include the Boolean semiring, the tropical
	semiring and various important generalisations thereof.
	Section~\ref{sec_greens}
	explores a consequence of exactness, generalising a result of \cite{K_tropd}
	to show that there is a straightforward characterisation of
	Green's $\GreenD$ relation for matrices with entries in an exact semiring.
	Finally, Section~\ref{sec_remarks} briefly discusses some questions arising
	from the results in the preceding sections.
	\section{Preliminaries}
	\label{sec_prelim}
	A \textit{semiring} is a commutative semigroup $(S,+)$ with an associative
	(but not necessarily commutative) multiplication $S\times S\to S$ that distributes
	over addition from both sides. In the literature it is generally assumed, see \cite{Golan99} for instance,
	that a semiring has an additive identity $0\in S$ that is an absorbing element
	for multiplication and a multiplicative identity $1\in S$.
	However we shall only require our semirings to satisfy the following (weaker)
	condition: for every non-empty finite subset $L\subseteq S$ there are \textit{local identities}
	$0_L,1_L\in S$ with $a+0_Lb=1_La=a$ and $a+b0_L=a1_L=a$ for all $a,b,\in L$.
	This condition is an adaptation of the `local zeroes' condition in \cite{K_tropd}.
	It allows us to consider (semirings derived from) the finitary tropical semiring
	(see Section~\ref{sec_anti}), which is of considerable importance in tropical
	algebraic geometry.

	We write $S^{m\times n}$ for the additive semigroup of $m$ row, $n$ column
	matrices with entries in $S$, where $m,n\in\mathbb N$. Matrix
	multiplication behaves in the usual ways: where defined it is associative and
	distributes over matrix addition. In particular $S^{n\times n}$ is a semiring,
	called the \textit{matrix semiring} $\mx_n(S)$, for each $n\in\mathbb N$.
	Note that local identities in $\mx_n(S)$ can be defined in terms of
	local identities in $S$, since the set of entries of a finite set
	of matrices is a finite subset of $S$.
	These local identities are analogues of the standard zero and identity matrices.

	Each matrix $A\in S^{m\times n}$ has an associated \textit{row space}
	\begin{align}
		\row(A)&=\bigl\{x\in S^{1\times n}:\text{$x=uA$ for some $u\in S^{1\times m}$}\bigr\}\\
		\intertext{and an associated \textit{column space}}
		\col(A)&=\bigl\{y\in S^{m\times 1}:\text{$y=Av$ for some $v\in S^{n\times 1}$}\bigr\}.
	\end{align}
	The local identities in $S$ ensure that $\row(A)$ actually contains the rows of $A$,
	and similarly that $\col(A)$ contains the columns of $A$. Notice that the row (column) space
	of $A$ is closed under addition and left (right) multiplication by $1\times 1$ matrices
	(which are nothing but elements of $S$). This suggests that the sets $\row(A)$
	and $\col(A)$ have `module' structures under the action of $S$.
	We now make this statement precise.

	A \textit{left $S$-module} $X$ is a commutative semigroup $(X,+)$ with an
	associative left action $S\times X\to X$ that distributes over addition in $X$ and $S$,
	subject to the requirement that for every non-empty finite $K\subseteq X$ there
	are left local identities $0_K,1_K\in S$ with $x+0_Ky=1_Kx=x$ for all $x,y\in K$.
	Right $S$-modules
	are similarly defined.
	In \cite{Golan99} S-modules are called `semimodules' and our
	local identity conditions are replaced by analogous conditions for the (global) identities in $S$.
	If $S$ is a (unital)
	ring then the present definition describes $S$-modules in the usual sense.
	In particular, if $S$ is a field then an $S$-module is just a vector space.

	It is clear that $S$ itself is both a left and a right $S$-module.
	We consider $S^{1\times n}$ and row spaces of matrices to be left $S$-modules,
	with the obvious left action, and we consider $S^{m\times 1}$ and
	column spaces of matrices to be right $S$-modules, again with the obvious right
	action.

	A function $\phi\colon X\to Y$ between left $S$-modules $X$ and $Y$ is called
	\textit{left linear} if $\phi(ax+by)=a(\phi x)+b(\phi y)$
	for all $x,y\in X$ and all $a,b\in S$. If a left linear (respectively
	\textit{right linear}; defined in the obvious way) function is bijective then
	the inverse function is automatically left (right) linear. Such a function $X\to Y$
	is called an \textit{isomorphism} and we write $X\cong Y$.
	If a linear function is injective then we call it an \textit{embedding}.

	The set $\col(A)^*$ of all right linear functions from the column space of a matrix
	$A\in S^{m\times n}$ to $S$ is a left $S$-module. The sum of two functions
	$\phi,\psi\colon\col(A)\to S$ is given by $(\phi+\psi)y=\phi y+\psi y$
	for all $y\in \col(A)$ and the left action of $S$ on $\phi\colon\col(A)\to S$
	is given by $(a\phi)y=a(\phi y)$ for all $a\in S$ and all $y\in\col(A)$.
	Left local identities for a finite non-empty $L\subseteq\col(A)^*$ are
	given by $0_{L'}$ and $1_{L'}$ where $L'=\{\phi A:\phi\in L\}$, and
	where $\phi A\in S^{1\times n}$ is the result of applying $\phi$ to the
	columns of $A$. Similarly the set $\row(A)^*$ of left linear functions
	$\row(A)\to S$ is a right $S$-module.
	\section{Exact semirings}
	\label{sec_exact}
	In this section we give the definition of exactness for semirings and we
	characterise it in ways that are familiar from classical
	linear algebra over fields. In particular we show that exactness is equivalent
	to the property of linear functions on row and column spaces extending to linear functions
	on the appropriate containing modules.
	\begin{definition}
		\label{def_exact}
		A semiring $S$ is \textit{exact} if the following conditions hold for
		all $A\in S^{m\times n}$.
		\begin{enumerate}
			\item[(E1)]
			If $x\in S^{1\times n}\setminus\row(A)$ then there are
			$v,v'\in S^{n\times 1}$ with $Av=Av'$ but $xv\neq xv'$.
			\item[(E2)]
			If $y\in S^{m\times 1}\setminus\col(A)$ then there are
			$u,u'\in S^{1\times m}$ with $uA=u'A$ but $uy\neq u'y$.
		\end{enumerate}
	\end{definition}
	Note that if $S$ is commutative then (E1) and (E2) are equivalent because
	the transpose of any matrix product is then the reverse product of the transposes
	of the factor matrices.

	Definition~\ref{def_exact} motivates a further two definitions.
	The \textit{kernel} of a set $X\subseteq S^{1\times n}$ of row vectors is
	the \textit{right congruence} (equivalence relation that is compatible with
	addition and the right action of $S$)
	\begin{equation}
		\label{eq_kernel}
		\ker(X)=\bigl\{(v,v')\in S^{n\times 1}\times S^{n\times 1}:\text{$xv=xv'$ for all $x\in X$}\bigr\}
	\end{equation}
	on $S^{n\times 1}$.
	In the case that $X$ is the row space of a matrix $A\in S^{m\times n}$, the
	kernel of $X$ is simply the (set-theoretic) kernel of the surjective right
	linear function $S^{n\times 1}\to\col(A)$ given by $v\mapsto Av$.
	This observation makes it clear that $S^{n\times 1}/\ker\row(A)\cong\col(A)$
	as right $S$-modules. The kernel of a set of column vectors is a similarly
	defined left congruence on $S^{1\times m}$, and we have that $S^{1\times m}/\ker\col(A)\cong\row(A)$
	as left $S$-modules.

	From the form of \eqref{eq_kernel} it is clear that if $X\subseteq Y\subseteq S^{1\times n}$
	then $\ker(Y)\subseteq\ker(X)$. In particular if $A\in S^{m\times n}$ and
	$B\in S^{p\times n}$ with $\row(A)\subseteq\row(B)$ then $\ker\row(B)\subseteq\ker\row(A)$.
	That is, $\ker(-)$ is inclusion-reversing for row (and similarly column) spaces
	of matrices. The following theorem tells us that $S$ is exact if and only $\ker(-)$ is an
	order anti-embedding (with respect to inclusion) for row and column spaces.
	\begin{theorem}
		\label{thm_kernels}
		Let $S$ be a semiring.
		Then $S$ is exact if and only if the following conditions hold for all
		$A\in S^{m\times n}$.
		\begin{enumerate}
			\item[(F1)]
			If $B\in S^{p\times n}$ with $\ker\row(A)\subseteq\ker\row(B)$ then
			$\row(B)\subseteq\row(A)$.
			\item[(F2)]
			If $B\in S^{m\times q}$ with $\ker\col(A)\subseteq\ker\col(B)$ then
			$\col(B)\subseteq\col(A)$.
		\end{enumerate}
		\begin{proof}
			We show the equivalence of (E1) and (F1) for all $A\in S^{m\times n}$,
			the equivalence of (E2) and (F2) being dual.
			First suppose that (E1) holds for $A\in S^{m\times n}$. To show that
			(F1) holds for $A$ let $B\in S^{p\times n}$ with $\ker\row(A)\subseteq\ker\row(B)$
			and let $x\in\row(B)$. Then $xv=xv'$ for all $(v,v')\in\ker\row(A)$,
			and as such $x\in\row(A)$ by the contrapositive of (E1).
			Hence $\row(B)\subseteq\row(A)$.

			Now suppose that (F1) holds for
			$A\in S^{m\times n}$. To show that (E1) holds for $A$ let
			$x\in S^{1\times n}\setminus\row(A)$. The matrix
			\begin{equation}
				B=\begin{bmatrix}
					A\\x
				\end{bmatrix}\in S^{(m+1)\times n}
			\end{equation}
			satisfies $\row(A)\subseteq\row(B)$, so
			$\ker\row(B)\subseteq\ker\row(A)$ because $\ker(-)$ is inclusion-reversing
			for row spaces. By construction $\row(B)\nsubseteq\row(A)$,
			so $\ker\row(B)\subset\ker\row(A)$ by the contrapositive of (F1).
			Therefore there is some $(v,v')\in\ker\row(A)$
			with $Bv\neq Bv'$, and hence with $xv\neq xv'$.
			That is, (E1) holds for $A$.
		\end{proof}
	\end{theorem}
	Given a matrix $A\in S^{m\times n}$ and any $x\in\row(A)$ we can define a
	right linear function $\col(A)\to S$ by $Av\mapsto xv$. This function is
	well-defined because $x$ can be written as $uA$ for some $u\in S^{1\times m}$,
	so $x\mapsto(Av\mapsto xv)$ is a well-defined function $\row(A)\to\col(A)^*$.
	This (outer) function is injective,
	since if $xv=x'v$ for all $v\in S^{n\times 1}$ then we can use local identities to deduce that $x=x'$,
	and is left linear,
	so is an embedding of left $S$-modules.
	Similarly there is an embedding $\col(A)\to\row(A)^*$ of right $S$-modules
	given by $y\mapsto(uA\mapsto uy)$.
	The following theorem characterises exactness in terms of the surjectivity of
	these embeddings.
	\begin{theorem}
		\label{thm_surjective}
		Let $S$ be a semiring.
		Then $S$ is exact if and only if the following conditions hold for all
		$A\in S^{m\times n}$.
		\begin{enumerate}
			\item[(G1)]
			The left embedding $\row(A)\to\col(A)^*$ given by $x\mapsto(Av\mapsto xv)$
			is surjective.
			\item[(G2)]
			The right embedding $\col(A)\to\row(A)^*$ given by $y\mapsto(uA\mapsto uy)$
			is surjective.
		\end{enumerate}
		\begin{proof}
			First suppose that (E1) holds for $A\in S^{m\times n}$.
			Let $\phi\in\col(A)^*$ and suppose that $\phi A\notin\row(A)$,
			where $\phi A\in S^{1\times n}$ is the result of applying $\phi$ to
			the columns of $A$.
			Then by (E1) there are $v,v'\in S^{n\times 1}$ with
			$Av=Av'$ and $(\phi A)v\neq(\phi A)v'$.
			However, right linearity of $\phi$ gives $\phi(Av)=(\phi A)v$ and
			$\phi(Av')=(\phi A)v'$, so we have $(\phi A)v=(\phi A)v'$.
			This contradicts $(\phi A)v\neq(\phi A)v'$, and as such we must
			actually have $\phi A\in\row(A)$.
			Therefore $\phi$ is the image of $\phi A$ under the given function
			$\row(A)\to\col(A)^*$.
			Hence (G1) holds for $A$.

			Now suppose that (G1) holds for
			$A\in S^{m\times n}$. We show that (F1) holds for $A$.
			Let $B\in S^{p\times n}$
			with $\ker\row(A)\subseteq\ker\row(B)$ and let $x\in\row(B)$.
			The right linear function $\col(A)\to S$ given by $Av\mapsto xv$ is
			then well-defined: if $Av=Av'$ then $(v,v')\in\ker\row(A)\subseteq\ker\row(B)$,
			so $xv=uBv=uBv'=xv'$ where $x=uB$ for some $u\in S^{1\times p}$.
			By (G1) this function $\col(A)\to S$ is given by
			$Av\mapsto x'v$ for some $x'\in\row(A)$, and as such $xv=x'v$ for all
			$v\in S^{n\times 1}$.
			We can then use local identities to deduce that $x=x'\in\row(A)$,
			and therefore $\row(B)\subseteq\row(A)$.
			Hence (F1) holds for $A$.

			As above, a similar argument shows that (G2) is equivalent,
			via (F2), to (E2) for all $A\in S^{m\times n}$.
		\end{proof}
	\end{theorem}
	Theorem~\ref{thm_surjective} tells us that if a semiring $S$ is exact then
	$\row(A)\cong\col(A)^*$ as left $S$-modules and $\col(A)\cong\row(A)^*$ as right
	$S$-modules for all $A\in S^{m\times n}$. It also tells us that every right linear
	function $\col(A)\to S$ can be written as an inner product, and so in particular
	every right linear function $\col(A)\to S$ extends to a right linear function
	$S^{m\times 1}\to S$. Similarly every left linear function $\row(A)\to S$
	extends to a left linear function $S^{1\times n}\to S$.
	In fact, taken together these properties turn out to be yet another characterisation of exactness.
	\begin{theorem}
		\label{thm_extend}
		Let $S$ be a semiring.
		Then $S$ is exact if and only if the following conditions hold for all
		$A\in S^{m\times n}$.
		\begin{enumerate}
			\item[(H1)]
			Every right linear function $\col(A)\to S$ extends to a right linear
			function $S^{m\times 1}\to S$.
			\item[(H2)]
			Every left linear function $\row(A)\to S$ extends to a left linear
			function $S^{1\times n}\to S$.
		\end{enumerate}
		\begin{proof}
			We have already observed that (H1)
			follows from (G1), so suppose that
			(H1) holds for $A\in S^{m\times n}$
			and let $\phi\in\col(A)^*$. Now let $I\in\mx_m(S)$ be an $m\times m$
			identity matrix local to $A$, that is, with $IA=A$.
			Recall that $\phi A$ is the result of applying $\phi$ to the columns
			of $A$. Since $\phi$ extends to $S^{m\times 1}$ we may also apply $\phi$
			to the columns of $I$, giving $\phi I$. Right linearity of $\phi$
			then gives $\phi(IAv)=(\phi I)Av$ for all $Av\in\col(A)$.
			Therefore $\phi(Av)=(\phi I)Av$ for all $Av\in\col(A)$ because $IA=A$,
			and as such $\phi$ is the image of $(\phi I)A$
			under the function $\row(A)\to\col(A)^*$ given in (G1).
			Hence (G1) holds for $A$.

			Similarly (H2) and (G2) are equivalent for all $A\in S^{m\times n}$.
		\end{proof}
	\end{theorem}
	Recall that a semiring $S$ is \textit{self-injective} if every linear
	function $X\to S$, for $X$ an arbitrary $S$-module, factors through every
	embedding of $X$ into any $S$-module. It is clear from Theorem~\ref{thm_extend}
	that exactness of $S$ is nothing but self-injectivity restricted to the embeddings
	$\col(A)\subseteq S^{m\times 1}$ and $\row(A)\subseteq S^{1\times n}$.
	Thus every self-injective semiring is exact.
	\section{Orthogonal complements and exact rings}
	\label{sec_rings}
	Throughout this section $R$ will denote a \textit{ring} (with unity, but not
	necessarily commutative), which for our purposes may be thought of as a semiring
	which has a global zero and identity (denoted $0$ and $1$) and which forms
	a group under addition (the inverse of $a\in R$ being denoted $-a$).

	The \textit{orthogonal complement} of a set $X\subseteq R^{1\times n}$ is
	the right $R$-module
	\begin{equation}
		X^\perp=\bigl\{v\in R^{n\times 1}:\text{$xv=0$ for all $x\in X$}\bigr\}.
	\end{equation}
	Similarly the orthogonal complement of a set $Y\subseteq R^{m\times 1}$ is
	the left $R$-module
	\begin{equation}
		Y^\perp=\bigl\{u\in R^{1\times m}:\text{$uy=0$ for all $y\in Y$}\bigr\}.
	\end{equation}
	In the case $n=m=1$ the notation `$^\perp$' is ambiguous, but the
	appropriate definition of orthogonal complement will always be clear from
	context. If in this case $R$ is commutative then the two notions of
	orthogonal complement coincide and $X^\perp$ is just the \textit{annihilator}
	$\ann(X)$ of $X\subseteq R$.

	If $A\in R^{m\times n}$ then $\row(A)^\perp$ is the kernel of the surjective
	right linear function $R^{n\times 1}\to\col(A)$ given by $v\mapsto Av$.
	Therefore $R^{n\times 1}/\row(A)^\perp\cong\col(A)$ as right $R$-modules.
	Similarly $R^{1\times m}/\col(A)^\perp\cong\row(A)$ as left $R$-modules.
	As with the kernels of row and column spaces, taking orthogonal complements
	reverses inclusions.
	We also have $\row(A)\subseteq\row(A)^{\perp\perp}$
	and $\col(A)\subseteq\col(A)^{\perp\perp}$ for all $A\in R^{m\times n}$,
	with (by the following theorem) equality if and only if $R$ is exact.
	\begin{theorem}
		\label{thm_orthogonal}
		Let $R$ be a ring.
		Then $R$ is exact if and only if $\row(A)^{\perp\perp}=\row(A)$
		and $\col(A)^{\perp\perp}=\col(A)$ for all $A\in R^{m\times n}$.
		\begin{proof}
			Suppose that (E1) holds for $A\in R^{m\times n}$ and let
			$x\in\row(A)^{\perp\perp}$. Then $xv=xv'$ for all
			$(v,v')\in\ker\row(A)$ because $v-v'\in\row(A)^\perp$. Therefore
			$x\in\row(A)$ by the contrapositive of (E1), and as such
			$\row(A)^{\perp\perp}\subseteq\row(A)$. Hence
			$\row(A)^{\perp\perp}=\row(A)$ because, as noted above,
			we always have $\row(A)\subseteq\row(A)^{\perp\perp}$.

			Now suppose that $\row(A)^{\perp\perp}=\row(A)$ for
			$A\in R^{m\times n}$ and let $B\in R^{p\times n}$ with
			$\ker\row(A)\subseteq\ker\row(B)$. Then
			$\row(A)^\perp\subseteq\row(B)^\perp$, so
			\begin{equation}
				\row(B)\subseteq\row(B)^{\perp\perp}\subseteq\row(A)^{\perp\perp}=\row(A)
			\end{equation}
			because taking orthogonal complements reverses inclusions.
			Hence (F1) holds for $A$.

			A similar argument shows that (E2) implies $\col(A)^{\perp\perp}=\col(A)$,
			and that this in turn implies (F2), for all $A\in R^{m\times n}$.
		\end{proof}
	\end{theorem}
	Recall that a ring is called an \textit{integral domain} if it is commutative, it has no zero
	divisors and $0\neq1$. The following proposition tells us
	that every exact integral domain must be a field, so in particular the
	prototypical integral domain $\mathbb Z$ is not exact.
	\begin{proposition}
		\label{prop_exact_commutative}
		If $R$ is an exact commutative ring then every non-zero element of $R$
		is a zero divisor or a unit.
		\begin{proof}
			Let $a\in R\setminus\{0\}$. If $a$ is a zero divisor then we are done,
			so suppose that $ab\neq0$ for all $b\in R\setminus\{0\}$. Then the function
			$\phi\in(aR)^*$ given by $\phi(av)=v$ is well-defined.
			Hence $1=\phi a=ua$ for some $u\in R$ by (G1).
		\end{proof}
	\end{proposition}
	We recall that a \textit{principal ideal domain} is an integral domain in which every ideal
	is generated by a single element. All fields, $\mathbb Z$ and $K[t]$ for $K$
	a field are principal ideal domains, where $K[t]$ denotes the ring of
	polynomials in $t$ with coefficients from $K$. For the remainder of this
	section $P$ will be a principal ideal domain and $R$ will be the quotient of
	$P$ by (the ideal generated by) a fixed element $r\in P$.
	\begin{theorem}
		\label{thm_pid_row_col}
		If $R=P/rP$ for $P$ a principal ideal domain and $r\in P$
		then $\row(A)\cong\col(A)$ for all $A\in R^{m\times n}$.
		\begin{proof}
			Adding rows and columns of zeroes to $A$ does not (up to isomorphism
			at least) change its row space or column space, so we may assume
			that $A$ is square.
			By \cite[Theorem 10.4]{Coh1974} there are invertible matrices $M$
			and $N$ with entries in $R$ and with $MAN$ diagonal, so we have
			\begin{equation}
				\begin{split}
					ANM^{-\transpose}&=M^{-1}MANM^{-\transpose}\\
					&=M^{-1}(MAN)^\transpose M^{-\transpose}\\
					&=M^{-1}N^\transpose A^\transpose\\
					&=\bigl(ANM^{-\transpose}\bigr)^\transpose
				\end{split}
			\end{equation}
			because $MAN$ is symmetric, and as such $ANM^{-\transpose}$ is symmetric.

			The required isomorphism
			$\phi\colon\row(A)\to\col(A)$ is given by
			$\phi(uA)=\bigl(uANM^{-\transpose}\bigr)^\transpose$ for all
			$uA\in\row(A)$, where the image of $\phi$ is contained in $\col(A)$
			because $ANM^{-\transpose}$ is symmetric.
			The inverse of $\phi$
			is given by $\phi^{-1}(Av)=\bigl(N^{-\transpose}MAv\bigr)^\transpose$
			for all $Av\in\col(A)$.
		\end{proof}
	\end{theorem}
	Theorem~\ref{thm_pid_row_col} places no restriction on $r\in P$, so it
	applies to (for instance) $\mathbb Z$ and $\mathbb Z/n\mathbb Z$ for all
	non-zero $n\in\mathbb Z$. The following two results also apply to all $\mathbb Z/n\mathbb Z$ with
	$n\neq0$.
	They do not, however, apply to $\mathbb Z$.
	\begin{theorem}
		\label{thm_nonzero}
		Let $R=P/rP$ for $P$ a principal ideal domain and $r\in P\setminus\{0\}$.
		If $A\in R^{m\times n}$ then there is some
		$B\in R^{n\times n}$ with $\row(A)^\perp=\col(B)$ and
		$\col(B)^\perp=\row(A)$.
		\begin{proof}
			We will first show that the result holds for all $A\in R^{1\times 1}$.
			That is, we will show that for all $A\in R$ there is some $B\in R$ with
			$\ann(RA)=BR$ and $\ann(BR)=RA$.

			Write $A=a'+rP$ for some $a'\in P$ and take $a\in P$ to generate the
			ideal $Pa'+rP$. That is, take $a$ to be the greatest common divisor
			of $a'$ and $r$ in $P$. Then $r=ba$ for some $b\in P$ and $RA=R(a+rP)$. If
			$a=0$ or $b=0$ then $r=0$, so we must have $a\neq0$ and $b\neq0$.
			Now write $B=b+rP$. Then we have $BR\subseteq\ann(RA)$ because
			$ba\in rP$. If $c+rP\in\ann(RA)$ then $ca=dr=dba$ for some
			$d\in P$, so $c=db$ because $P$ is an integral domain and $a\neq0$.
			Therefore $c+rP\in BR$, and as such $\ann(RA)=BR$. Similarly
			$\ann(BR)=RA$.
			Note that if $A=0+rP$ then
			we can take $B=1+rP$.

			This result extends to diagonal matrices $A\in R^{m\times n}$.
			Without loss of generality (by adding or removing rows of zeroes as
			necessary) we may assume that $A$ is square diagonal.
			To form $B\in R^{n\times n}$ replace each diagonal entry $a\in R$
			of $A$ by an element $b\in R$ using the above procedure.
			This will ensure
			that $\row(A)^\perp=\col(B)$ and $\col(B)^\perp=\row(A)$.

			We can now prove the theorem for an arbitrary matrix
			$A\in R^{m\times n}$. By \cite[Theorem 10.4]{Coh1974}
			there are invertible matrices $M\in\mx_m(R)$
			and $N\in\mx_n(R)$ with $MAN$ diagonal, so by the above result
			there is a matrix
			$B\in R^{n\times n}$ with $\row(MAN)^\perp=\col(B)$ and
			$\col(B)^\perp=\row(MAN)$.
			Invertibility of $M$ and $N$ then gives
			\begin{gather}
				\row(A)^\perp=\row(MA)^\perp=\col(NB)\\
				\intertext{and}
				\col(NB)^\perp=\row(MA)=\row(A).
			\end{gather}
		\end{proof}
	\end{theorem}
	\begin{corollary}
		\label{cor_exact}
		If $R=P/rP$ for $P$ a principal ideal domain and $r\in P\setminus\{0\}$
		then $R$ is exact and
		\begin{equation}
			\row(A)^\perp\cong R^{1\times n}/\row(A)
		\end{equation}
		for all $A\in R^{m\times n}$.
		\begin{proof}
			By Theorem~\ref{thm_nonzero} there is some $B\in R^{n\times n}$
			with $\row(A)^\perp=\col(B)$ and $\col(B)^\perp=\row(A)$.
			Therefore $\row(A)^{\perp\perp}=\col(B)^\perp=\row(A)$.
			We could verify, using a result dual to Theorem~\ref{thm_nonzero},
			that $\col(A)^{\perp\perp}=\col(A)$ too, but there is no need since $R$
			is commutative.
			Hence $R$ is exact by Theorem~\ref{thm_orthogonal}.

			For the second claim, by Theorem~\ref{thm_pid_row_col} we have
			$\row(A)^\perp=\col(B)\cong\row(B)$.
			Therefore $\row(A)^\perp\cong R^{1\times n}/\col(B)^\perp=R^{1\times n}/\row(A)$.
		\end{proof}
	\end{corollary}
	If $r\neq0$ then $R$ is known to be self-injective \cite[Theorem~4.35]{Rot1979},
	so $R$ is exact by Theorem~\ref{thm_extend}. Nevertheless, Theorem~\ref{thm_nonzero}
	and Corollary~\ref{cor_exact} are included here as the proof of Theorem~\ref{thm_nonzero}
	describes a procedure to
	compute the orthogonal complement of the row space of a matrix with entries in $R$.
	\section{Matrix semirings and group semirings}
	\label{sec_matrix}
	A subset $T$ of a semiring $S$ will be called a \textit{subsemiring} of $S$
	if it is closed under addition and multiplication, and if for every non-empty
	finite $L\subseteq S$ there are $0_L,1_L\in T$ with $a+0_Lb=1_La=a$ and
	$a+b0_L=a1_L=a$ for all $a,b,\in L$.
	That is, for $T$ to be a subsemiring of $S$ we must be able to choose local
	identities that lie in $T$, not just in $S$
	(as the present definition of a semiring requires), for every non-empty finite subset of $S$.
	This is by analogy with the standard definition \cite{Golan99} of a subsemiring $T$, which
	requires $T$ to contain $0,1\in S$.

	It is clear from the relevant definitions that if $T$ is a subsemiring of $S$
	then $S$ is both a left and a right $T$-module.
	A subsemiring $T$ of $S$ will be called a \textit{left} (\textit{right})
	\textit{retract} of $S$ if there is a left (right) $T$-linear function
	$S\to T$ that fixes $T$ pointwise.
	\begin{theorem}
		\label{thm_retract}
		Let $S$ be an exact semiring and $T$ be a (right and left) retract of
		$\mx_n(S)$. If $S^{1\times n}$ and $S^{n\times 1}$ embed into $T$,
		as right and left $T$-modules respectively, then $T$ is exact.
		\begin{proof}
			To apply Theorem~\ref{thm_kernels} (which will establish the
			exactness of $T$)
			we need to verify that (F1) and (F2)
			hold for all $A\in T^{m\times q}$.
			We will only show that (F1) holds for all $A$, the proof that (F2) holds
			being dual.

			Let $A\in T^{m\times q}$.
			Since $T\subseteq\mx_n(S)$ we can view $A$ as a matrix with entries
			in $S$ that is divided into blocks of size $n\times n$.
			That is, we can view $A$ as an element of $S^{mn\times qn}$.
			We write $\row_S(A)$ for the row space
			of $A$ when considered as a matrix with entries in $S$
			and to prevent ambiguity we write $\row_T(A)$ for the row space of
			$A\in T^{m\times q}$.

			Let $B\in T^{p\times q}$ and suppose that
			$\ker\row_T(A)\subseteq\ker\row_T(B)$.
			We will show that $\ker\row_S(A)\subseteq\ker\row_S(B)$.
			By assumption there is an injective left $T$-linear function
			$\phi\colon S^{n\times 1}\to T$. Left $T$-linearity of $\phi$ gives
			\begin{equation}
				A(\phi v)=
				\begin{bmatrix}
					\phi(X_1v)\\
					\vdots\\
					\phi(X_mv)
				\end{bmatrix}
			\end{equation}
			for all $v\in S^{qn\times 1}$, where $X_1,\dotsc,X_m\in T^{1\times q}$
			are the rows of $A$ and $\phi v\in T^{q\times 1}$ is the result of
			applying $\phi$ to the $q$ blocks of $v$.
			Therefore if $v,v'\in S^{qn\times 1}$ with $Av=Av'$ then
			$A(\phi v)=A(\phi v')$
			because $X_iv=X_iv'$ for each $1\leq i\leq m$.
			That is, if $(v,v')\in\ker\row_S(A)$ then
			$(\phi v,\phi v')\in\ker\row_T(A)$.
			In fact, as $\phi$ is injective we have
			\begin{equation}
				(v,v')\in\ker\row_S(A)\quad\Leftrightarrow\quad
				(\phi v,\phi v')\in\ker\row_T(A)
			\end{equation}
			for all $v,v'\in S^{qn\times 1}$, and similarly for $B$.
			Hence $\ker\row_S(A)\subseteq\ker\row_S(B)$
			because $\ker\row_T(A)\subseteq\ker\row_T(B)$.

			Since $S$ is exact, (F1) holds for $A\in S^{mn\times qn}$.
			Therefore $\row_S(B)\subseteq\row_S(A)$ because $\ker\row_S(A)\subseteq\ker\row_S(B)$.
			To verify that (F1) holds for $A\in T^{m\times q}$ we need to show that
			$\row_T(B)\subseteq\row_T(A)$, so let $X\in\row_T(B)\subseteq T^{1\times q}$.
			Then each row $x\in S^{1\times qn}$ of $X$ is in $\row_S(B)\subseteq\row_S(A)$,
			so each row $x\in S^{1\times qn}$ of $X$ can be written as $x=uA$ for
			some $u\in S^{1\times mn}$.
			Therefore $X=UA$ for some $U\in S^{n\times mn}$.

			By assumption again there is a right $T$-linear function $\pi\colon\mx_n(S)\to T$
			that fixes $T$ pointwise. Right $T$-linearity of $\pi$ gives
			\begin{equation}
				(\pi U)A=
				\begin{bmatrix}
					\pi(UY_1)&\cdots&\pi(UY_q)
				\end{bmatrix},
			\end{equation}
			where $Y_1,\dotsc,Y_q\in T^{m\times 1}$ are the columns of $A$
			and $\pi U\in T^{1\times m}$ is the result of applying $\pi$ to the $m$ blocks of $U$.
			Since $UA=X\in T^{1\times q}$, the $q$ entries of $X$ are
			$UY_1,\dotsc,UY_q\in T$. Therefore
			\begin{equation}
				X=
				\begin{bmatrix}
					UY_1&\cdots&UY_q
				\end{bmatrix}
				=(\pi U)A
			\end{equation}
			because $\pi$ fixes $T$ pointwise, and as such $X\in\row_T(A)$.
			Hence $\row_T(B)\subseteq\row_T(A)$, as required for (F1) to hold
			for $A\in T^{m\times q}$.
		\end{proof}
	\end{theorem}
	Theorem~\ref{thm_retract} gives sufficient conditions for a
	(sub)semiring of matrices to be exact. In particular it tells us that if
	$S$ is exact then every matrix semiring $\mx_n(S)$ will be be exact too, because $S^{1\times n}$ and
	$S^{n\times1}$ clearly embed into $\mx_n(S)$.
	\begin{corollary}
		If $S$ is an exact semiring then every matrix semiring $\mx_n(S)$ is exact.
	\end{corollary}
	Now let $S$ be a semiring and $(\Gamma,\cdot,1)$ be a finite group.
	The set of functions $\Gamma\to S$ is a semiring with operations
	defined by $(f+g)\alpha=f\alpha+g\alpha$ and
	\begin{equation}
		\label{group_mult}
		(fg)\alpha=\sum_{\beta\gamma=\alpha}f\beta\cdot g\gamma
		=\sum_{\beta\in\Gamma}f\beta\cdot g\bigl(\beta^{-1}\alpha\bigr)
	\end{equation}
	for all $f,g\colon\Gamma\to S$ and all $\alpha\in\Gamma$.
	We call this semiring the \textit{group semiring} $S\Gamma$.
	Local identities for a finite non-empty subset $L\subseteq S\Gamma$
	can be defined in terms of $0_{L'},1_{L'}\in S$, where
	$L'=\{f\alpha:\text{$f\in L$ and $\alpha\in\Gamma$}\}$.
	Notice that $L'$ is finite because $L$ and $\Gamma$ are finite.
	\begin{theorem}
		\label{thm_group_exact}
		If $S$ is an exact semiring and $\Gamma$ is a finite group then the
		group semiring $S\Gamma$ is exact.
		\begin{proof}
			Let $\Gamma$ be a finite group with $\lvert\Gamma\rvert=n$ and
			$\Gamma=\{\alpha_1,\dotsc,\alpha_n\}$.
			We begin by showing that $S\Gamma$ can be identified
			with a subsemiring of $\mx_n(S)$ via the map sending a function
			$f\in S\Gamma$
			to the matrix
			\begin{equation}
				\label{eq_group_matrix}
				F=
				\begin{bmatrix}
					f\bigl(\alpha_1^{-1}\alpha_1\bigr)&\cdots&f\bigl(\alpha_1^{-1}\alpha_n\bigr)\\
					\vdots&\ddots&\vdots\\
					f\bigl(\alpha_n^{-1}\alpha_1\bigr)&\cdots&f\bigl(\alpha_n^{-1}\alpha_n\bigr)
				\end{bmatrix}.
			\end{equation}
			Notice that we can recover $f$ from any row or column of $F$,
			so the map sending $f$ to $F$ is injective.
			If $f,g\in S\Gamma$ with corresponding matrices $F,G\in\mx_n(S)$ then
			\begin{equation}
				F+G=
				\begin{bmatrix}
					(f+g)\bigl(\alpha_1^{-1}\alpha_1\bigr)&\cdots&(f+g)\bigl(\alpha_1^{-1}\alpha_n\bigr)\\
					\vdots&\ddots&\vdots\\
					(f+g)\bigl(\alpha_n^{-1}\alpha_1\bigr)&\cdots&(f+g)\bigl(\alpha_n^{-1}\alpha_n\bigr)
				\end{bmatrix}
			\end{equation}
			by the definition of addition in $S\Gamma$, so $F+G$ is the matrix
			corresponding to $f+g\in S\Gamma$.
			It therefore remains to show that $FG$ is the matrix corresponding
			to $fg\in S\Gamma$.
			For each $1\leq i,j\leq n$ we have
			\begin{equation}
				(FG)_{ij}=\sum_{k=1}^nF_{ik}\cdot G_{kj}
				=\sum_{k=1}^nf\bigl(\alpha_i^{-1}\alpha_k\bigr)\cdot g\bigl(\alpha_k^{-1}\alpha_j\bigr)
			\end{equation}
			by \eqref{eq_group_matrix}, so
			\begin{equation}
				(FG)_{ij}=\sum_{\beta\in\Gamma}f\beta\cdot g\bigl(\beta^{-1}\alpha_i^{-1}\alpha_j\bigr)
				=(fg)\bigl(\alpha_i^{-1}\alpha_j\bigr)
			\end{equation}
			for each $1\leq i,j\leq n$ by \eqref{group_mult}.
			Hence $FG$ is the matrix corresponding to $fg$.
			We can therefore treat $S\Gamma$ as if it were the subsemiring of $n\times n$
			matrices of the form \eqref{eq_group_matrix}.

			Note that $S\Gamma$ is a subsemiring of $\mx_n(S)$ in the sense defined above
			because we can take local identities for any subset of matrices
			that are analogues of the standard zero and identity matrices.
			Therefore we can take local identities that are of the form \eqref{eq_group_matrix}.

			We will ultimately use Theorem~\ref{thm_retract} to show that $S\Gamma\subseteq\mx_n(S)$
			is exact, but first we introduce some notation.
			Let $\rho\colon\mx_n(S)\to S^{1\times n}$ and $\kappa\colon\mx_n(S)\to S^{n\times 1}$
			be the functions that select the first row and first column respectively
			of a matrix $A\in\mx_n(S)$.
			It is clear that $\rho$ is right $\mx_n(S)$-linear, so in particular
			it is right $S\Gamma$-linear.
			Similarly $\kappa$ is left $S\Gamma$-linear.
			Now let $\psi\colon S^{1\times n}\to S\Gamma$ be the function that
			sends $x\in S^{1\times n}$ to the unique, by \eqref{eq_group_matrix},
			element of $S\Gamma$ with first row $x$.
			Similarly let $\phi\colon S^{n\times 1}\to S\Gamma$ be the function
			that sends $v\in S^{n\times 1}$
			to the unique element of $S\Gamma$ with first column $v$.

			If $x\in S^{1\times n}$ then $\psi x$ has first row $x$, so $\rho\psi x=x$.
			On the other hand, if $F\in S\Gamma$ then $\psi\rho F=F$ because
			$\psi\rho F$ and $F$ are in $S\Gamma$ and have the same first row.
			Therefore $\rho$ restricted to $S\Gamma$ has inverse $\psi$,
			and as such $\psi$ is left $S\Gamma$-linear because $\rho$ is.
			Similarly $\kappa$ restricted to $S\Gamma$ has inverse $\phi$,
			so $\phi$ is right $S\Gamma$-linear.

			We are now ready to apply Theorem~\ref{thm_retract}.
			We first show that $S\Gamma\subseteq\mx_n(S)$ is a right
			and a left retract of $\mx_n(S)$. For the right $S\Gamma$-linear
			function $\mx_n(S)\to S\Gamma$ take $\psi\circ\rho$ and for the
			left $S\Gamma$-linear function $\mx_n(S)\to S\Gamma$ take $\phi\circ\kappa$.
			As noted above, these functions fix $S\Gamma$ pointwise, so $S\Gamma$
			is both a right and a left retract of $\mx_n(S)$.
			For the right $S\Gamma$-module embedding $S^{1\times n}\to S\Gamma$
			take $\psi$ (which is injective by definition) and for the left $S\Gamma$-module
			embedding $S^{n\times 1}\to S\Gamma$ take $\phi$.
			Theorem~\ref{thm_retract} then ensures that $S\Gamma$ is exact.
		\end{proof}
	\end{theorem}
	\section{Tropical, Boolean and other anti-involutive semirings}
	\label{sec_anti}
	In this section we define what it means for a semiring to be anti-involutive
	and we show that all such semirings are exact.
	We also show that the row and column spaces of matrices with entries in an
	anti-involutive semiring are, in a sense defined below, anti-isomorphic as
	modules.

	Recall that a semiring $S$ is \textit{idempotent} if $a+a=a$ for all $a\in S$.
	If in this case $X$ is an $S$-module then it follows that $x+x=x$ for all $x\in X$ too,
	and as such we can define a partial order on $X$ by setting
	\begin{equation}
		x\leq y\quad\Leftrightarrow\quad x+y=y
	\end{equation}
	for all $x,y\in X$.
	In particular $S$ itself can be partially ordered in this way.
	Like addition, the partial order on $S$ extends entrywise to matrices
	(of the same size) with entries in $S$.
	\begin{example}
		\label{ex_tropical}
		Let $(G,\cdot,e)$ be a (torsion-free) group with a total order that is compatible with the group operation,
		by which we mean that if $a,b\in G$ with $a\leq b$
		then $ca\leq cb$ and $ac\leq bc$ for all $c\in G$.
		For instance, $G$ could be $(\mathbf Z,+,0)$, $(\mathbf Q,+,0)$ or $(\mathbf R,+,0)$ with the standard ordering.
		We can turn $G$ into an idempotent semiring by defining
		\begin{equation}
			a+b=\max\{a,b\}
		\end{equation}
		for all $a,b\in G$.
		Notice that the partial order induced by this operation is just the original total order on $G$.
		Local identities for a non-empty finite $L\subseteq G$
		are $1_L=e$ and
		\begin{equation}
			0_L=\min\bigl\{(\min L)(\max L)^{-1},(\max L)^{-1}(\min L)\bigr\}.
		\end{equation}
		In order to distinguish it from the group $G$ we write $G_{\mathrm{max}}$ for this semiring.
	\end{example}
	A function $\phi\colon X\to Y$ between left $S$-modules $X$ and $Y$,
	for $S$ an idempotent semiring, will be called \textit{left monotone} if
	\begin{equation}
		ax\leq x'\quad\Rightarrow\quad a(\phi x)\leq\phi x'
	\end{equation}
	for all $x,x'\in X$ and all $a\in S$.
	It is clear that if $\phi$ is left linear (e.g., matrix multiplication on the left)
	then $\phi$ is left monotone,
	but the converse is not true in general.
	For example, the function $\phi\colon\mathbb R_{\mathrm{max}}^{1\times 2}\to\mathbb R_{\mathrm{max}}$
	given by
	\begin{equation}
		\phi\begin{bmatrix}
			x_1&x_2
		\end{bmatrix}
		=-\max\{-x_1,-x_2\}
	\end{equation}
	is left monotone but not left linear since
	\begin{equation}
		1=
		\phi\max\Bigl\{
		\begin{bmatrix}
			1&0
		\end{bmatrix},
		\begin{bmatrix}
			0&1
		\end{bmatrix}
		\Bigr\}
		\neq
		\max\Bigl\{
		\phi\begin{bmatrix}
			1&0
		\end{bmatrix},
		\phi\begin{bmatrix}
			0&1
		\end{bmatrix}
		\Bigr\}=0.
	\end{equation}
	However, if $\phi$ is left monotone
	and has a left monotone inverse then $\phi$ (and hence also $\phi^{-1}$) is
	left linear, as the following proposition shows.
	\begin{proposition}
		\label{prop_isomorphism}
		Let $S$ be an idempotent semiring and $\phi\colon X\to Y$ be a function
		between left (right) $S$-modules.
		Then $\phi$ is an isomorphism of left (right) $S$-modules if and only if
		$\phi$ is left (right) monotone and has a left (right) monotone inverse.
		\begin{proof}
			If $\phi$ is an isomorphism then
			$\phi$ and $\phi^{-1}$ are left linear, so are left monotone.
			It therefore remains to show the converse.
			Let $x,y\in X$ and $a,b\in S$.
			Then $ax\leq ax+by$ and $by\leq ax+by$, so
			$a(\phi x)\leq\phi(ax+by)$ and $b(\phi y)\leq\phi(ax+by)$
			because $\phi$ is left monotone.
			Therefore $a(\phi x)+b(\phi y)\leq\phi(ax+by)$.
			We also have $a(\phi x)\leq a(\phi x)+b(\phi y)$ and
			$b(\phi y)\leq a(\phi x)+b(\phi y)$, so
			\begin{align}
				ax&\leq\phi^{-1}\bigl(a(\phi x)+b(\phi y)\bigr)\\
				\intertext{and}
				by&\leq\phi^{-1}\bigl(a(\phi x)+b(\phi y)\bigr)
			\end{align}
			because $\phi^{-1}$ is left monotone.
			Therefore $\phi(ax+by)\leq a(\phi x)+b(\phi y)$,
			because $\phi$ is order-preserving,
			with $a(\phi x)+b(\phi y)\leq\phi(ax+by)$ from above.
			Hence $\phi$ is left linear, and as such $\phi$ is an isomorphism
			of left $S$-modules.
		\end{proof}
	\end{proposition}
	In response to Proposition~\ref{prop_isomorphism},
	a function $\phi\colon X\to Y$ from a left $S$-module $X$ to a right
	$S$-module $Y$ will be called \textit{antitone} if
	\begin{equation}
		\label{antitone}
		ax\leq x'\quad\Rightarrow\quad(\phi x')a\leq\phi x
	\end{equation}
	for all $x,x'\in X$ and all $a\in S$.
	Antitone functions are (in particular) order-reversing.
	If $\phi$ is antitone and has an antitone inverse $\phi^{-1}\colon Y\to X$
	then we call $\phi$ (and also $\phi^{-1}$) an \textit{anti-isomorphism}
	and we say that $X$ and $Y$ are \textit{anti-isomorphic} (as modules). The condition for
	$\phi^{-1}$ to be antitone is just the obvious analogue of \eqref{antitone}:
	\begin{equation}
		\label{antitone_dual}
		ya\leq y'\quad\Rightarrow\quad a\bigl(\phi^{-1}y'\bigr)\leq\phi^{-1}y
	\end{equation}
	for all $y,y'\in Y$ and all $a\in S$.
	Note that by Proposition~\ref{prop_isomorphism} the composition of two anti-isomorphisms is an isomorphism of
	left (or right, whichever is appropriate) modules.

	For the remainder of this section $S$ will be an idempotent semiring with an
	involutive anti-isomorphism $\overline{\phantom-}\colon S\to S$. We call such
	semirings \textit{anti-involutive}.
	Ordinarily, specifying that a function $\phi\colon S\to S$ is an anti-isomorphism
	is ambiguous if $S$ is not commutative: should \eqref{antitone} and \eqref{antitone_dual} apply as written,
	or should $\phi$ and $\phi^{-1}$ be interchanged?
	However, the fact that $\overline{\phantom-}$ is an involution removes this
	distinction, meaning that
	$\overline{\phantom-}$ is antitone in the sense of \eqref{antitone} and
	in the sense of \eqref{antitone_dual}.

	The semirings $G_{\mathrm{max}}$ in Example~\ref{ex_tropical} are anti-involutive
	with $\overline a=a^{-1}$ for all $a\in G$.
	In particular the \textit{finitary tropical} semiring $\ft=\mathbb R_{\mathrm{max}}$
	is anti-involutive with $\overline a=-a$ for all $a\in\mathbb R$.
	Another anti-involutive semiring is the \textit{Boolean} semiring $\mathbb B=\{\bot,\top\}$,
	where addition is given by `or', multiplication by `and' and the involution by `not'.

	The involution on $S$ extends to matrices with entries in $S$.
	Specifically, given $A\in S^{m\times n}$ we define $\overline A\in S^{n\times m}$ by
	$\overline A_{ij}=\overline{A_{ji}}$ for all $1\leq i\leq m$ and all $1\leq j\leq n$.
	Proposition~\ref{prop_cycle} (below) tells us that $\overline{\phantom-}$ is
	an involutive anti-isomorphism on $\mx_n(S)$
	for all $n$, so every matrix semiring $\mx_n(S)$ is anti-involutive.
	Every group semiring $S\Gamma$ for $\Gamma$ a finite group is anti-involutive
	with $\overline f\alpha=\overline{f\alpha^{-1}}$ for $f\in S\Gamma$
	and all $\alpha\in\Gamma$.
	This can also be seen by identifying $S\Gamma$ with a subsemiring of matrices
	(as in the proof of Theorem~\ref{thm_group_exact}) that is closed under the
	involution just defined.

	In fact, if $S$ is \textit{complete} in the sense of order theory
	(so that sums of arbitrarily many elements are possible, such as in $\mathbb B$)
	then Proposition~\ref{prop_cycle} actually holds for infinite matrices and, by implication, elements
	of $S\Gamma$ for an infinite group $\Gamma$.
	We will not discuss this idea further here, except to say that the following example can
	be applied to an arbitrary group $\Gamma$.
	\begin{example}
		A group semiring $\mathbb B\Gamma$ can be viewed as the semiring of subsets of $\Gamma$,
		where addition is union and multiplication is given by $AB=\{\alpha\beta:\text{$\alpha\in A$ and $\beta\in B$}\}$
		for all $A,B\subseteq\Gamma$.
		Such a semiring is anti-involutive with $\overline A=\Gamma\setminus A^{-1}$
		for all $A\subseteq\Gamma$, where $A^{-1}=\bigl\{\alpha^{-1}:\alpha\in A\bigr\}$.
	\end{example}
	\begin{proposition}
		\label{prop_cycle}
		If $S$ is an anti-involutive semiring then
		\begin{equation}
			MA\leq B\quad\Rightarrow\quad\overline BM\leq\overline A
		\end{equation}
		for all $A\in S^{m\times n}$, all $B\in S^{p\times n}$ and all $M\in S^{p\times m}$.
		\begin{proof}
			Suppose that $MA\leq B$.
			We will show that $\bigl(\overline BM\bigr)_{ji}\leq\overline{A}_{ji}$
			for all $1\leq i\leq m$ and all $1\leq j\leq n$.
			Since $MA\leq B$ we have $M_{ki}A_{ij}\leq(MA)_{kj}\leq B_{kj}$
			for all $1\leq k\leq p$, so $\overline{B_{kj}}M_{ki}\leq\overline{A_{ij}}$
			for all $k$ because $\overline{\phantom-}\colon S\to S$
			is antitone.
			Therefore $\overline B_{jk}M_{ki}\leq\overline A_{ji}$ for all $k$,
			and as such
			\begin{equation}
				\bigl(\overline BM\bigr)_{ji}=\sum_{k=1}^p\overline B_{jk}M_{ki}\leq\overline A_{ji}.
			\end{equation}
			Hence $\overline BM\leq\overline A$.
		\end{proof}
	\end{proposition}
	An immediate consequence of Proposition~\ref{prop_cycle} is that
	\begin{equation}
		\label{matrix_antitone}
		MA\leq B\quad\Rightarrow\quad
		\overline BM\leq\overline A\quad\Rightarrow\quad
		A\overline B\leq\overline M\quad\Rightarrow\quad
		MA\leq B
	\end{equation}
	for all $A\in S^{m\times n}$, all $B\in S^{p\times n}$ and all $M\in S^{p\times m}$.
	In the arguments that follow we make frequent use of \eqref{matrix_antitone}.
	Informally, it allows us to cycle terms in an inequality at the cost of
	introducing a $\overline{\phantom-}$ to the two terms that cross the inequality.
	We also use the fact that matrix multiplication is order-preserving.
	\begin{theorem}
		\label{thm_anti-isomorphic}
		If $S$ is an anti-involutive semiring then $\row(A)$ and $\col(A)$ are
		anti-isomorphic for all $A\in S^{m\times n}$.
		\begin{proof}
			We exhibit an antitone function $\phi\colon\row(A)\to\col(A)$
			that has an antitone inverse $\col(A)\to\row(A)$. Define $\phi$ by
			$\phi x=A\overline x$ for all $x\in\row(A)$. If $ax\leq x'$ for
			$x,x'\in\row(A)$ and $a\in S$ then $\overline{x'}a\leq\overline x$,
			and as such $(\phi x')a\leq\phi x$. Hence $\phi$ is antitone.
			Similarly the function $\psi\colon\col(A)\to\row(A)$ given by
			$y\mapsto\overline yA$ is antitone, so it remains to show that
			$\psi$ is the inverse of $\phi$.
			Let $uA\in\row(A)$. Then
			\begin{equation}
				\psi\phi(uA)=\overline{A\overline{uA}}A\leq uA
			\end{equation}
			by \eqref{matrix_antitone} because $A\overline{uA}\leq A\overline{uA}$.
			We also have
			$uA\leq uA$, so by \eqref{matrix_antitone} again $A\overline{uA}\leq\overline u$,
			and thus $u\leq\overline{A\overline{uA}}$.
			Therefore
			\begin{equation}
				uA\leq\overline{A\overline{uA}}A=\psi\phi(uA),
			\end{equation}
			and as such $\psi\phi(uA)=uA$. Similarly $\phi\psi(Av)=Av$ for all
			$Av\in\col(A)$, so $\psi$ is the inverse of $\phi$.
		\end{proof}
	\end{theorem}
	In \cite{K_tropd} it is proved that the row and column spaces of matrices
	with entries in $\ft$ are anti-isomorphic.
	The notion of anti-isomorphism used in \cite{K_tropd} is, in general, weaker
	than the present one, but in the case of $\ft$ the two are equivalent.
	In \cite{Cohen04} it is only proved that if $\mathcal K$ is a complete idempotent
	semiring with an anti-isomorphism $\mathcal K\to\mathcal K$ then the row and
	column spaces of any matrix with entries in $\mathcal K$ are anti-isomorphic
	as complete lattices, and not (as in Theorem~\ref{thm_anti-isomorphic}) as modules.
	This is because the anti-isomorphism $\mathcal K\to\mathcal K$ need not be
	an involution.

	The completion of $\ft$ is obtained from $\ft$ by adjoining
	two new elements $-\infty$ and $\infty$ with $-\infty<a<\infty$ for all $a\in\ft$,
	setting $(-\infty)a=a(-\infty)=-\infty$ for all $a\in\ft\cup\{-\infty,\infty\}$
	and $(\infty)a=a(\infty)=\infty$ for all $a\in\ft\cup\{\infty\}$.
	The row and column spaces of matrices with entries in this semiring are known to be
	anti-isomorphic in the respective senses above \cite{Cohen04,K_tropd},
	but we conclude that they are in fact anti-isomorphic in the strongest sense because
	the completion of $\ft$ is anti-involutive (with $\overline{-\infty}=\infty$ and $\overline{\infty}=-\infty$).
	\begin{theorem}
		If $S$ is an anti-involutive semiring then $S$ is exact.
		\begin{proof}
			Let $A\in S^{m\times n}$ and $B\in S^{p\times n}$ with
			$\ker\row(A)\subseteq\ker\row(B)$. We will show that
			$\row(B)\subseteq\row(A)$.
			Let $uB\in\row(B)$ and for convenience
			let $v=\overline{uB}$. Then $\overline{Av}A\leq\overline v=uB$
			by \eqref{matrix_antitone} because $Av\leq Av$.
			As in the proof of Theorem~\ref{thm_anti-isomorphic},
			\eqref{matrix_antitone} gives $A\overline{\overline{Av}A}=Av$, so
			\begin{equation}
				\Bigl(\overline{\overline{Av}A},v\Bigr)\in\ker\row(A)\subseteq\ker\row(B).
			\end{equation}
			Therefore $B\overline{\overline{Av}A}=Bv\leq\overline u$, and as such
			$uB\leq\overline{Av}A$ by \eqref{matrix_antitone} again.
			Hence $uB=\overline{Av}A\in\row(A)$, as required
			for $\row(B)\subseteq\row(A)$. Therefore (F1) and, by a similar argument, (F2)
			hold for all $A\in S^{m\times n}$. Thus $S$ is exact by
			Theorem~\ref{thm_kernels}.
		\end{proof}
	\end{theorem}
	\section{Green's relations}
	\label{sec_greens}
	The equivalence relations of Green are usually used to reveal
	the structure of semigroups \cite{Green51,Howie95}, but they can also be
	defined for (not necessarily square) matrices.
	Two matrices $A\in S^{m\times n}$ and $B\in S^{p\times n}$ are
	$\GreenL$-related, written $A\GreenL B$, if there are $M\in S^{p\times m}$
	and $P\in S^{m\times p}$ with $MA=B$ and $PB=A$. We notice that
	$A\GreenL B$ if and only if $\row(A)=\row(B)$. Similarly
	$A\in S^{m\times n}$ and $B\in S^{m\times q}$ are $\GreenR$-related, written
	$A\GreenR B$, if there are $N\in S^{n\times q}$ and $Q\in S^{q\times n}$
	with $AN=B$ and $BQ=A$. Again, we notice that $A\GreenR B$ if and only if
	$\col(A)=\col(B)$.

	Green's $\GreenD$ relation is the relational composition of $\GreenL$ and $\GreenR$.
	That is, $A\in S^{m\times n}$ and $B\in S^{p\times q}$ are $\GreenD$-related
	if there is some $C\in S^{p\times n}$ with $A\GreenL C\GreenR B$.
	\begin{theorem}
		\label{thm_green}
		Let $S$ be an exact semiring, $A\in S^{m\times n}$ and $B\in S^{p\times q}$.
		Then the following conditions are equivalent.
		\begin{enumerate}
			\item
			\label{thm_green_row}
			$\row(A)\cong\row(B)$ as left $S$-modules.
			\item
			\label{thm_green_d}
			$A\GreenD B$.
			\item
			\label{thm_green_col}
			$\col(A)\cong\col(B)$ as right $S$-modules.
		\end{enumerate}
		\begin{proof}
			We will only show the equivalence of conditions (\ref{thm_green_d}) and
			(\ref{thm_green_col}), the equivalence of (\ref{thm_green_row}) and
			(\ref{thm_green_d}) being dual.

			If $A\GreenD B$ then there is some $C\in S^{p\times n}$ with $A\GreenL C\GreenR B$.
			Therefore there are $M\in S^{p\times m}$ and $P\in S^{m\times p}$ with $MA=C$ and $PC=A$.
			There are also $N\in S^{n\times q}$ and $Q\in S^{q\times n}$ with $CN=B$ and $BQ=C$.
			The right linear function $\col(A)\to\col(B)$ given by $Av\mapsto MAv=BQv$
			has inverse $\col(B)\to\col(A)$ given by $Bv\mapsto PBv=ANv$ because $PMA=A$ and $MPB=B$.
			Hence $\col(A)\cong\col(B)$.

			Now let $\phi\colon\col(A)\to\col(B)$ be an isomorphism of right $S$-modules.
			For each $1\leq i\leq p$ the composition of $\phi$ with projection to the $i$th entry
			yields a right linear function $\phi_i\in\col(A)^*$.
			By (G1), Theorem~\ref{thm_surjective},
			each $\phi_i$ is given by $Av\mapsto u_iAv$ for some $u_i\in S^{1\times m}$.
			Therefore $\phi$ is given by $Av\mapsto MAv$ with $\col(MA)=\col(B)$, where
			\begin{equation}
				M=
				\begin{bmatrix}
					u_1\\
					\vdots\\
					u_p
				\end{bmatrix}\in S^{p\times m}.
			\end{equation}
			A similar argument shows that $\phi^{-1}$ is given by $Bv\mapsto PBv$
			for some $P\in S^{m\times p}$ with $PMA=A$.
			Therefore $A\GreenL MA\GreenR B$, and as such $A\GreenD B$.
		\end{proof}
	\end{theorem}
	Two matrices $A\in S^{m\times n}$ and $B\in S^{p\times q}$ are $\GreenJ$-related
	if and only if there are matrices $M\in S^{p\times m}$, $N\in S^{n\times q}$ with
	$MAN=B$ and matrices $P\in S^{m\times p}$, $Q\in S^{q\times n}$ with $PBQ=A$.
	It is clear from the definitions that ${\GreenD}\subseteq{\GreenJ}$, but even
	if $S$ is exact it need not be the case that ${\GreenD}={\GreenJ}$.
	For instance, it is known that ${\GreenD}={\GreenJ}$ for matrices with entries
	in $\ft$ and that ${\GreenD}\neq{\GreenJ}$ for matrices with entries in its
	completion \cite{K_tropj}.
	It is also known that ${\GreenD}={\GreenJ}$ for matrices with entries in any
	Artinian ring that has a maximal ideal \cite{Cao2010}.
	\section{Remarks and open questions}
	\label{sec_remarks}
	In Section~\ref{sec_rings} we showed that proper quotients of principal ideal
	domains are exact rings. As remarked there, this result also follows from our
	characterisation of exactness as a weak form of self-injectivity
	(Theorem~\ref{thm_extend}), since proper
	quotients of principal ideal domains are self-injective rings.
	In fact, all known (to the authors at least) examples of exact rings are also
	self-injective.
	This leaves open the question of whether there exist any exact rings that are
	not self-injective.

	Corollary~\ref{cor_exact} not only establishes the exactness of proper quotients
	of principal ideal domains, it also provides a description of the orthogonal complement
	of the row space of a matrix.
	However, as with the relationship between row and column spaces (isomorphism in Section~\ref{sec_rings} vs.\ anti-isomorphism in Section~\ref{sec_anti}),
	it seems unlikely that a general description of the kernels of matrices is
	obtainable using exactness alone.
	A natural extension of the present work would therefore be to determine the
	structure of the kernels of matrices with entries in an anti-involutive semiring.
	\bibliographystyle{plain}
	\bibliography{references}
	\end{document}